\documentclass{amsart}

\usepackage{amssymb,latexsym}
\usepackage{palatino}
\usepackage{mathrsfs}

\newcommand{\GL}{\operatorname{GL}}
\newcommand{\SL}{\operatorname{SL}}

\newcommand{\Hom}{\operatorname{Hom}}
\newcommand{\End}{\operatorname{End}}
\newcommand{\Aut}{\operatorname{Aut}}
\newcommand{\tensor}{\otimes}
\newcommand{\soc}{\operatorname{soc}}

\newcommand{\Dist}{\operatorname{Dist}}
\newcommand{\Aff}{\mathbf{A}}

\newcommand\Z{\mathbf{Z}}    
\newcommand\F{\mathbf{F}}    
\newcommand\G{\mathbf{G}}    

\DeclareMathAlphabet{\mathrc}{U}{eur}{m}{n}
\newcommand\Witt{\mathrc{W}}   

\newcommand\Lie{\operatorname{Lie}}
\newcommand\Ad{\operatorname{Ad}}

\newcommand\lie[1]{\mathfrak{#1}}

\newcommand\glie{\lie{g}}

\newcommand\mm{\mathfrak{m}}

\newcommand{\normal}{\lhd}
\newcommand{\semidirect}{\rtimes}

\newcommand{\congruent}{\equiv} 
\newcommand{\iso}{\simeq}

\newtheorem{theorem}{Theorem}
\newtheorem{prop}[theorem]{Proposition}
\newtheorem{cor}[theorem]{Corollary}
\newtheorem{lem}[theorem]{Lemma}
\theoremstyle{remark}
\newtheorem{rem}[theorem]{Remark}
\newtheorem{rems}[theorem]{Remarks}
\newtheorem{example}[theorem]{Example}

\newcommand{\A}{\mathscr{A}}

\newcommand{\Int}{{\operatorname{Int}}}

\newcommand{\XX}{\mathscr{X}}
\newcommand{\YY}{\mathscr{Y}}
\newcommand{\ZZ}{\mathscr{Z}}

\begin{document}
\bibliographystyle{hamsalpha}

\title[Faithful representations]{Faithful representations
of $\SL_2$ over truncated Witt vectors} 
\author{George J. McNinch}
\thanks{This work was supported in part by a grant from the National
  Science Foundation.}

\address{Department of Mathematics \\
         Room 255 Hurley Building \\
         University of Notre Dame \\
         Notre Dame, Indiana 46556-5683 \\
         USA}
\email{mcninch.1@nd.edu}
\date{February 11, 2002}

\begin{abstract}
  Let $\Gamma_2$ be the six dimensional linear algebraic $k$-group
  $\SL_2(\Witt_2)$, where $\Witt_2$ is the ring of Witt vectors of
  length two over the algebraically closed field $k$ of characteristic
  $p>2$.  Then the minimal dimension of a faithful rational
  $k$-representation of $\Gamma_2$ is $p+3$.
\end{abstract}

\maketitle
 
\section{Introduction}

Let $\Witt = \Witt(k)$ be the ring of Witt vectors over the
algebraically closed field $k$ of characteristic $p>0$. Let $\Witt_n =
 \Witt/p^n\Witt$ be the ring of length $n$ Witt
vectors. (See \cite[II.\S6]{serre79:_local_field} for definitions and
basic properties of Witt vectors, and see \S \ref{sec:witt-vectors}
below.)  We regard $\Witt_n$ as a ``ring variety'' over $k$, the
underlying variety of which is $\Aff^n_{/k}$.  If $n \ge 2$, the ring
$\Witt_n$ is not a $k$-algebra.

Let $\Gamma_n = \SL_2(\Witt_n)$ be the group of $2 \times 2$ matrices
with entries in $\Witt_n$ and with determinant 1.  Then $\Gamma_n$ is
a closed subvariety of the $4n$ dimensional affine space of $2 \times
2$ matrices over $\Witt_n$; thus, it is an affine algebraic group over
$k$.  As such, it is a closed subgroup of $\GL(V)$ for some finite
dimensional $k$-vector space $V$, i.e. it has a faithful finite
dimensional $k$-linear representation. Note that for $n \ge 2$,
$\Witt_n$ is not a vector space over $k$ in any natural way, so the
natural action of $\Gamma_n$ on $\Witt_n \oplus
\Witt_n$ is not a $k$-linear representation.

Let $H$ be any linear algebraic group over $k$.  A rational $H$-module
$(\rho,V)$ is said to be faithful if $\rho$ defines a closed embedding
$H \to \GL(V)$; this is equivalent to the condition: both $\rho$ and
$d\rho$ are injective.

\begin{theorem}
  \label{theorem:min-degree}
  If $(\rho,V)$ is a representation of $\Gamma_2$ with $\dim V \le
  p+2$, then $\rho(u^p) = 1_V$ for each unipotent element $u \in
  \Gamma_2$.
\end{theorem}
\begin{theorem}
  \label{theorem:min-construction}
  If $p>2$, 
  the minimal dimension of a faithful rational representation of
  $\Gamma_2$ is $p+3$. 
\end{theorem}

With the same notation, if $p=2$ then $\Gamma_2$ has a rational
representation $(\rho,V)$ with $\dim V = p+3 = 5$, and $\rho$ is
abstractly faithful (i.e. injective on the closed points of $G$) but
$\ker d\rho$ is the Lie algebra of a maximal torus of $\Gamma_2$.

After some preliminaries in \S 2 through \S 4, we construct in \S
\ref{sec:faithful-example} a representation $(\rho,V)$ of $\Gamma_2$
of dimension $p+3$ and show that $\rho$ is abstractly faithful; in \S
\ref{sec:inf-faithful} we show finally that $d\rho$ is injective when
$p>2$. Combined with Theorem \ref{theorem:min-degree}, this proves
Theorem \ref{theorem:min-construction}.

In \S \ref{sec:minimality} we prove that the unipotent radical $R$ of
$\Gamma_2$ acts trivially on any rational module with dimension $\le
p+2$; this completes the verification of Theorem
\ref{theorem:min-degree}. An important tool in the proof is a result
obtained in \S \ref{sec:distinct-weight-spaces} concerning the weight
spaces of a representation of the group $\Witt_2 \semidirect
k^\times$; this result is proved with the help of the algebra of
distributions of the unipotent group $\Witt_2$.

Finally, in \S \ref{sec:finite-case}, we prove the analogue of Theorem
\ref{theorem:min-construction} for the finite groups $\Gamma_2(\F_q)$
provided that $p>2$ and $q \ge p^2$. In its outline, the proof is the
same as in the algebraic case. In the finite case, we replaced the
arguments concerning the algebra of distributions of $\Witt_2$ in \S
\ref{sec:distinct-weight-spaces} with some more elementary arguments
(see Proposition \ref{prop:finite-weight-spaces}). In fact, we could
use these more elementary arguments in the ``algebraic'' case, but the
techniques in \S \ref{sec:distinct-weight-spaces} give more
information and are therefore perhaps of independent interest. Note that
the condition on $q$ is an artifact of the proof; I do not know
if $\Gamma_2(\F_p) = \SL_2(\Z/p^2\Z)$ has a faithful $k$-representation
of dimension $< p+3$.

Thanks to Jens Carsten Jantzen and Jean-Pierre Serre for some helpful
comments on this manuscript.

\section{A negative application: unipotent elements in reductive groups.}

Let $H$ be a connected reductive group over $k$, and let $u \in H$ be
unipotent of order $p$.  If $p$ is a good prime for $H$,
there is a homomorphism $\SL_2(k) \to
H$ with $u$ in its image. This was proved by
Testerman \cite{testerman}; see also \cite{Mc:sub-principal}.

Now suppose that $u$ is a unipotent element in $H$ with order $p^n$,
$n \ge 1$.  Then there is a homomorphism $\Witt_n = \G_a(\Witt_n) \to
H$ with $u$ in its image.  This was proved by Proud \cite{Proud-Witt};
see \cite{mcninch-math.RT/0007056} for another proof when $H$ is classical.

In view of these results, one might wonder whether $u$ lies in the
image of a homomorphism $\gamma:\SL_2(\Witt_n) \to H$.  Theorem
\ref{theorem:min-degree} shows that, in general, the answer is ``no''.

Indeed, let $H$ be the reductive group $\GL_{p+1/k}$. Then a regular
unipotent element $u$ of $H$ has order $p^2$. On the other hand, if
$f:\SL_2(\Witt_2) \to H$ is a homomorphism, the theorem shows that $u$
is not in the image of $f$.

\section{Witt vectors}
\label{sec:witt-vectors}

Elements of $\Witt_n$ will be represented as tuples
$(a_0,a_1,\ldots,a_{n-1})$ with $a_i \in k$. For $w=(a_0,a_1)$ and
$w'=(b_0,b_1)$ in $\Witt_2$, we have:
\begin{equation}
  \label{eq:witt-vector-ops}
  w + w' = (a_0 + b_0,a_1+b_1 + F(a_0,b_0)) \quad \text{and} \quad
  w\cdot w' = (a_0b_0,a_0^pb_1 + b_0^pa_1),
\end{equation}
where $F(X,Y) = (X^p + Y^p - (X+Y)^p)/p \in \Z[X,Y]$.  

We have also the identity in $\Witt_n$
\begin{equation}
  \label{eq:torus-action-on-W_n}
  (t,0,\cdots,0) \cdot (a_0,a_1,\dots,a_{n-1}) = 
  (ta_0,t^pa_1,\dots,t^{p^{n-1}}a_{n-1})
\end{equation}
for all $t \in k$ and $(a_0,\cdots,a_{n-1}) \in \Witt_n$.

Let $\XX_n:\Witt_n = \G_a(\Witt_n) \to \Gamma_n$ and $\phi:k^\times \to \Gamma_n$ be the
maps
\begin{equation*}
\XX_n(w) = 
\begin{pmatrix}
  1 & w \\
  0 & 1
\end{pmatrix} \quad \text{and} \quad \phi(t) = 
\begin{pmatrix}
  (t,0,\dots,0) & 0 \\
  0 & (1/t,0,\dots,0)
\end{pmatrix}. 
\end{equation*}
Using \eqref{eq:torus-action-on-W_n}, one observes the relation
\begin{equation}
  \label{eq:conj}
  \Int(\phi(t))\XX_n(a_0,a_1,\dots,a_{n-1}) = 
  \XX_n(t^2a_0,t^{2p}a_1,\dots,t^{2p^{n-1}}a_{n-1}).
\end{equation}
The element $\XX_n(a_0,a_1,\dots,a_{n-1})$ is unipotent; if $a_0 \not
= 0$, it has order $p^n$.

For $n \ge 2$, the map $(a_0,a_1,\dots,a_{n-1}) \mapsto
(a_0,\dots,a_{n-2}):\Witt_n \to \Witt_{n-1}$ induces a surjective
homomorphism
\begin{equation*}
  \Gamma_n \to \Gamma_{n-1}
\end{equation*}
whose kernel we denote by $R_n$. Similarly, the
residue map $(a_0,a_1,\dots,a_{n-1}) \mapsto a_0:\Witt_n \to k$ induces
a surjective homomorphism $\eta:\Gamma_n \to \SL_2(k)$.

Concerning $R_n$, we have the following:
\begin{lem}
  \label{lem:adjoint-radical}
  The group $R_n$ is a connected, Abelian unipotent group of dimension
  3. More precisely, there is a $\Gamma_n$-equivariant isomorphism of algebraic
  groups
  \begin{equation}
    \label{eq:adjoint-radical}
    \gamma:R_n \to \lie{sl}_2(k);
  \end{equation}
  the action of $\Gamma_n$ on $\lie{sl}_2(k)$ is by $\Ad^{[n-1]} \circ
  \eta$, where $\Ad^{[n-1]}$ is the $(n-1)$-st Frobenius twist of the
  adjoint representation of $\SL_2(k)$, and the action of $\Gamma_n$ on
  $R_n$ is by inner automorphisms.
\end{lem}

The lemma follows from \cite[II.\S4.3]{Demazure-Gabriel}. Actually, the
cited result is quite straightforward for $\SL_2$.

The lemma shows that the kernel of $\eta$ is a $3(n-1)$ dimensional
unipotent group. In particular, $\Gamma_n$ has dimension $3n$. Since
$\Gamma_n/\ker \eta \iso \SL_2(k)$ is reductive, $\ker \eta$ is the
unipotent radical of $\Gamma_n$.  In particular, we see that the image of
$\phi$ is a maximal torus $T$ of $\Gamma_n$.

We now consider the case $n=2$; we write $R$ for $R_2$.
Let $\ZZ:\G_a(k) \to R < G$ be the homomorphism 
\begin{equation}
  \label{eq:ZZ-defined}
  \ZZ(s) = 
  \begin{pmatrix}
    (1,s) & 0 \\
    0 & (1,-s) 
  \end{pmatrix}.
\end{equation}
An easy matrix calculation yields:
\begin{equation}
  \label{eq:centralize}
  \Int(\phi(t))\ZZ(s) = \ZZ(s) \quad \text{for each}\ t \in k^\times 
  \ \text{and}\ s \in k.
\end{equation}

Recall that any non-0 nilpotent element of $\lie{sl}_2(k)$ is a cyclic
generator as a $\Gamma_2$-module, and any non-0 semisimple element of
$\lie{sl}_2(k)$ generates the socle of this module. (These remarks are
trivial for $p>2$ since in that case $\lie{sl}_2(k)$ is a simple
$\SL_2(k)$-module; the assertions in characteristic 2 are well known
and anyhow easy to verify). We thus obtain the following:
\begin{lem}
  \label{lem:generate}
  There are no proper $\Gamma_2$-invariant subgroups of $R$ containing
  $\XX_2(0,1)$. Any non-trivial $\Gamma_2$-invariant subgroup of $R$
  contains $\ZZ(1)$.
\end{lem}

\begin{rem}
  \label{rem:no-Levi}
  J. Humphreys pointed out that $\Gamma_n$ provides an example of a linear
  algebraic group in characteristic $p$ with no Levi decomposition.
  Here is an argument for his observation using the main result of
  this paper.
  
  First, since $\Gamma_2$ is a quotient of $\Gamma_n$, the above observation
  follows from:
  \begin{itemize}
  \item   The group $\Gamma_2$ has no Levi decomposition.
  \end{itemize}
  
  Let $H = (\Ad^{[1]},\lie{sl}_2(k)) \semidirect \SL_2(k)$. If we know
  that $H$ has a representation $(\mu,V)$ such that $\dim V < p+3$
  and $\ker \mu$ is finite, then Theorem \ref{theorem:min-degree}
  implies that $\Gamma_2$ is not isomorphic to $H$, hence that $\Gamma_2$ has no
  Levi decomposition.
  
  If $(\lambda,V)$ is a rational representation of a linear algebraic
  group $A$, we may form the semidirect product $\hat A = (\lambda,V)
  \semidirect A$.  There is a rational representation $(\hat \lambda,
  V \oplus k)$ of $\hat A$ given by $\hat \lambda(v,a)
  (w,\alpha) = (\lambda(a)w + \alpha v,\alpha)$ for $(v,a) \in \hat A$
  and $(w,\alpha) \in V \oplus k$. A straightforward check yields
  \begin{equation*}
    \ker \hat \lambda = \{ (0,a) \mid a \in \ker \lambda\}.
  \end{equation*}
  
  Applying this construction with $A = \SL_2(k)$, $(\lambda,V) =
  (\Ad^{[1]},\lie{sl}_2(k))$, $\hat A = H$, we find a representation
  \begin{equation*}
    (\widehat{\Ad^{[1]}},\lie{sl}_2(k) \oplus k)
  \end{equation*}
  with dimension $4<p+3$ and finite kernel $\{(0,\pm 1)\} \le H$, as
  required.

  For a different proof of this observation (when $p \ge 5$) see
  \cite[\S IV.23]{serre-l-adic}.
\end{rem}

\begin{rem}
  \label{rem:normal-subgroups}
  One can list all normal subgroups $N$ of $\Gamma_2$. If $p>2$, $N \cap R$
  must be either $1$ or $R$. Since the only proper, non-trivial normal
  subgroup of $\SL_2(k) = \Gamma_2/R$ is $\{\pm 1\}$, we see that $N$ is
  one of
  \begin{equation*}
    \Gamma_2, \quad R, \quad \{\pm 1\}, \quad R \cdot \{\pm 1\}, \quad 1
  \end{equation*}
  In this case $\{\pm 1\}$ is the center of $\Gamma_2$.
  
  If $p=2$, $N \cap R$ is either $1$, $R$, or $Z$, the inverse image
  under $\gamma$ of the 1 dimensional center of $\lie{sl}_2(k)$.  The
  group $\SL_2(k) = \Gamma_2/R$ is (abstractly) a simple group.  Let $N
  \normal \Gamma_2$ satisfy $N \cap R = Z$. Then $\eta(N)$ is either
  trivial or equal to $\SL_2(k)$.  If $\eta(N) \not = 1$, there is an
  extension
  \begin{equation*}
    1 \to Z \to N \to \SL_2(k) \to 1.
  \end{equation*}
  We have $H^2(\SL_2(k),Z) = H^2(\SL_2(k),k) = 0$ by \cite[Proposition
  II.4.13]{JRAG}, so such an extension must be split. But a splitting
  would yield a Levi decomposition for $\Gamma_2$, contrary to our
  observations in Remark \ref{rem:no-Levi}. Thus $\eta(N)=1$ so $N=Z$.
  [Note that the argument we just gave depends on Theorem
  \ref{theorem:min-degree}; we will not use it in proving this
  theorem.]

  To summarize, the possibilities for $N$ are:
  \begin{equation*}
    \Gamma_2, \quad R, \quad Z, \quad 1
  \end{equation*}
  The group $Z$ is the center of $G$. It is equal to the image
  $\ZZ(k)$.
\end{rem}

\section{Unipotent radicals and representations}
\label{sec:unip-radical}

\newcommand{\RR}{\mathcal{R}}

Let $A$ be a linear algebraic group over $k$, and let $R$ denote its
unipotent radical. If $(\rho,V)$ is a rational finite dimensional
$A$-representation (with $V \not = 0$), then the space $V^R$ of
$R$-fixed points is a non-0 $A$-subrepresentation (the fact that it is
non-0 follows from the Lie-Kolchin Theorem \cite[Theorem
6.3.1]{springer98:_linear_algeb_group}). This implies that there is a
filtration of $V$ by $A$-subrepresentations
\begin{equation}
  \label{eq:uni-radical-filtration}
  V=\RR^0V \supset \RR^1V \supset \RR^2 V \supset \cdots
  \supset \RR^n V = 0   
\end{equation}
with the properties: $(\rho(x)-1)\RR^i V \subset \RR^{i+1} V$ for each
$x \in R$ and each $i$, and each quotient $\RR^i V / \RR^{i+1} V$ is a
non-0 representation for the reductive group $A/R$.  

We see in particular that the simple $A$-modules are precisely the
simple $A/R$-modules inflated to $A$.

All this applies especially for $A= \Gamma_n$, $n \ge 1$.  We identify the
simple $\SL_2(k)$ modules and the simple $\Gamma_n$-modules; for $a \ge 0$,
there is thus a simple $\Gamma_n$-module $L(a)$ with highest weight $a$. If
$0 \le a \le p-1$, $\dim L(a) = a+1$. If $a$ has $p$-adic expansion $a
= \sum a_i p^i$ where $0 \le a_i \le p-1$ for each $i$, then
Steinberg's tensor product theorem \cite[II.3.17]{JRAG} yields
\begin{equation*}
  L(a) \iso L(a_0) \tensor L(a_1)^{[1]} \tensor L(a_1)^{[2]} \tensor \cdots 
\end{equation*}
where $V^{[i]}$ denotes the $i$-th Frobenius twist of the $\Gamma_n$ module
$V$.

\section{A faithful $G$-representation}
\label{sec:faithful-example}

In this section, we consider the group $G = \Gamma_2 = \SL_2(\Witt_2)$.  We
recall the homomorphisms $\XX_2:\Witt_2 \to G$ and $\ZZ:k \to R$; we
write $\XX$ for $\XX_2$.

\begin{lem}
  \label{lem:faithful}
  Let $(\rho,V)$ be a rational finite dimensional $G$-representation.
  Then $\rho$ is abstractly faithful (i.e. injective on the closed
  points of $G$) if and only if (i) $(\rho_{\mid T},V)$ is an
  abstractly faithful $T$-representation, and (ii) $u = \rho(\ZZ(1))
  \not = 1_V$.
\end{lem}

\begin{proof}
  The necessity of conditions (i) and (ii) is clear, so suppose these
  conditions hold and let $K$ be the kernel of $\rho$.  Let
  $\operatorname{gr}(V)$ denote the associated graded space for any
  filtration as in \eqref{eq:uni-radical-filtration}. Then
  $\operatorname{gr}(V)$ is a module for $\Gamma_2/R=\SL_2(k)$.  Condition
  (i) implies that $\operatorname{gr}(V)$ is an abstractly faithful representation
  of $\SL_2(k)$. Thus $\ker \rho$ is contained in $R$.  Now (ii)
  together with Lemma \ref{lem:generate} imply that $\ker \rho = 1$ as
  desired.
\end{proof}

\begin{theorem}
  \label{theorem:faithful}
  $G$ has an abstractly faithful representation $(\rho,V)$ with $\dim
  V = p+3$.
\end{theorem}

\begin{proof}
  Since $G$ acts as a group of automorphisms on the 4 dimensional
  affine $k$-variety $\Witt_2 \oplus \Witt_2$, there is a
  representation $\rho$ of $G$ on the coordinate ring $\A=k[\Witt_2
  \oplus \Witt_2]$ given by $(\rho(g)f)(w) = f(g^{-1}w)$ for $g \in G,
  f \in \A, w \in \Witt_2 \oplus \Witt_2$.  If we denote by $A_0$ and
  $A_1$ the coordinate functions on $\Witt_2 \oplus 0$, and by $B_0$
  and $B_1$ those on $0 \oplus \Witt_2$, then $\A$ identifies with the
  polynomial ring $k[A_0,A_1,B_0,B_1]$.
  
  There is a linear representation $\lambda$ of $k^\times$ on $\A$
  given by $(\lambda(t)f)(w) = f((t,0).w)$ for $t \in k^\times, f \in
  \A, w \in \Witt_2 \oplus \Witt_2$. One checks easily that
  $\lambda(t)A_0 = tA_0$, and that $\lambda(t)A_1 = t^pA_1$ for $t \in
  k^\times$, with similar statements for $B_0$ and $B_1$.
  
  For $\nu \in \Z$, let $\A_\nu$ be the space of all functions $f \in
  \A$ for which $\lambda(t)f = t^\nu f$ for all $t \in k^\times$ (i.e.
  the $\nu$-weight space for the torus action $\lambda$). Then we have
  a decomposition $\A= \bigoplus_{\nu \in \Z} \A_\nu$ as a
  $\lambda(k^\times)$-representation.
  
  Since $G$ acts ``$\Witt_2$-linearly'' on $\Witt_2 \oplus \Witt_2$,
  $\lambda(k^\times)$ centralizes $\rho(G)$; thus each $\A_\nu$ is a
  $G$-subrepresentation of $\A$. We consider the $G$-representation
  $(\rho_p,\A_p)$. One sees that $\A_p$ is spanned by all $A_0^i
  B_0^j$ with $i+j=p$ and $i,j \ge 0$, together with $A_1$ and $B_1$.
  Thus $\dim \A_p = p + 3$.
  
  Using \eqref{eq:witt-vector-ops}, one checks for each $s \in k$ that
  \begin{equation*}
    \rho_p(\ZZ(s)) A_1 = A_1 + sA_0^p,  
  \end{equation*}
  so that $\rho_p(\ZZ(1)) \neq 1$.  Since $A_1$ has $T$-weight $p$,
  $\A_p$ is an abstractly faithful representation of $T$; thus the lemma shows
  that $(\rho_p,\A_p)$ is an abstractly faithful $G$-representation.
\end{proof}

\begin{rems}
  \label{rem:faithful-remarks}
  \begin{itemize}
  \item[(a)] It is straightforward to see that $(\rho_p,\A_p)$ has
    length three, and that its composition factors are $L(p-2)$
    together with two copies of $L(p) = L(1)^{[1]}$.   
  
    
  \item[(b)] The representation $\rho_p$ is defined over the prime
    field $\F_p$.  In particular, the finite group $\SL_2(\Z/p^2\Z)$
    has a faithful representation on a $p+3$ dimensional $\F_p$-vector
    space.  More generally, the finite group $\SL_2(\Witt_2(\F_q))$
    has a faithful representation on a $p+3$ dimensional $\F_q$-vector
    space for each $q = p^r$).
    
  \item[(c)] We will show in \S \ref{sec:inf-faithful} that the representation
    $(\rho_p,\A_p)$ is actually faithful provided that $p>2$.
  \end{itemize}
\end{rems}

\section{Algebras of distributions}
\label{sec:distinct-weight-spaces}

Let $H$ be a linear algebraic $k$-group, and let $\Dist(H)$ be the
algebra of distributions on $H$ supported at the identity; see
\cite[I.7]{JRAG} for the definitions. Recall that elements of
$\Dist(H)$ are certain linear forms on the coordinate algebra $k[H]$.

The algebra structure of $\Dist(H)$ is determined by the
comultiplication $\Delta$ of $k[H]$; the product of $\mu,\nu \in
\Dist(H)$ is given by
\begin{equation*}
  \mu\cdot\nu:k[H] \xrightarrow{\Delta} k[H] \tensor_k k[H] 
  \xrightarrow{\mu \tensor \nu} k \tensor_k k = k.
\end{equation*}
We immediately see the following:
\begin{equation}
  \text{If $H$ is Abelian, then $\Dist(H)$ is a commutative $k$-algebra}.
\end{equation}

Now consider the case $H = \Witt_2$. As a variety, $\Witt_2$
identifies with $\Aff^2_{/k}$. We write $k[\Witt_2] = k[A_0,A_1]$ as
before.  As a vector space $\Dist(\Witt_2)$ has a basis
$\{\gamma_{i,j} \mid i,j \ge 0\}$ where $\gamma_{i,j}(A_0^sA_1^t) =
\delta_{i,s}\delta_{j,t}$; see \cite[I.7.3]{JRAG}.

Let $(\rho,V)$ be a $\Witt_2$-representation. This is determined by a
comodule map
\begin{equation*}
  \Delta_V:V \to V \tensor_k k[\Witt];
\end{equation*}
for $v \in V$ we have $\Delta_V(v) = \sum_{i,j \ge 0} \psi_{i,j}(v)
\tensor A_0^iA_1^j$ where $\psi_{i,j} \in \End_k(V)$.

The $\Witt_2$-module $(\rho,V)$ becomes a $\Dist(\Witt_2)$-module by
the recipe give in \cite[I.7.11]{JRAG}. A look at that recipe shows
that the basis elements $\gamma_{i,j} \in \Dist(\Witt_2)$ act on $V$
as multiplication by $\psi_{i,j}$. Since $\Witt_2$ is Abelian, we
deduce that the linear maps $\{\psi_{i,j} \mid i,j \ge 0\}$ pairwise
commute.

In view of the commutativity, we obtain
\begin{equation*}
  1_V=\rho(a,b)^{p^2}=\left(\sum_{i,j \ge 0} a^ib^j \psi_{i,j}\right)^{p^2} =
  \sum_{i,j \ge 0} a^{ip^2}b^{jp^2} \psi_{i,j}^{p^2}
\end{equation*}
identically in $a,b$; thus $\psi_{0,0} = 1_V$ and $\psi_{i,j}^{p^2} = 0$
if $i > 0$ or $j>0$.

Now let $H$ be the subgroup of $G = \SL_2(\Witt_2)$ generated by the
maximal torus $T$ together with the image of $\XX_2:\Witt_2 \to G$. Thus
$H$ is a semidirect product $\XX_2(\Witt_2) \semidirect T$.

Let $(\rho,V)$ be an $H$-representation.  The $T$-module structure on
$V$ yields a $T$-module structure on $\End_k(V)$; for a weight $\mu$
of $T$ we have $\psi \in \End_k(V)_\mu$ if and only if $\psi(v) \in
V_{\lambda + \mu}$ for all weights $\lambda$ and all $v \in
V_\lambda$.

Fix a weight vector $v \in V_\lambda$. Then
\begin{equation*}
  \rho(\XX_2(a,b))v = \sum_{i,j \ge 0} a^ib^j\psi_{i,j}(v),
\end{equation*}
where the $\psi_{i,j}$ are determined as before by the comodule map
for the $\Witt_2$-module $V$. A look at \eqref{eq:conj} shows that
$\psi_{i,j}(v) \in V_{\lambda + 2i + 2pj}.$ It follows that
$\psi_{i,j} \in \End_k(V)_{2i + 2pj}$.

\begin{prop}
  \label{prop:weight-estimates}
  Let $(\rho,V)$ be a representation of $H = \XX_2(\Witt_2)
  \semidirect T$.  Suppose that $\rho(\XX_2(0,1)) \not = 1_V$. Then
  $T$ has at least $p+1$ distinct weights on $V$. More precisely,
  there are weights $s \in \Z_{>0}$ and $\lambda \in \Z$ such that
  $V_{\lambda + 2sj} \not = 0$ for $0 \le j \le p$.
\end{prop}

\begin{proof}
  We have $\XX_2(0,1) = \XX_2(1,0)^p$. With notation as above, our
  hypothesis means that
  \begin{equation*}
    1_V \not = \rho(\XX_2(1,0))^p = \left(\sum_{i \ge 0} \psi_{i,0}\right)^p = 
    \sum_{i \ge 0} \psi_{i,0}^p.    
  \end{equation*}
  Thus there is some $s >0$ for which $\psi_{s,0}^p \not = 0$. Write
  $\psi = \psi_{s,0}$.  Recall that $\psi$ has $T$-weight $2s$. We may
  find a weight $\lambda \in \Z$ and $v \in V_\lambda$ for which
  $\psi^p(v) \not = 0$.  But then $v,\psi(v),\dots,\psi^p(v)$ are all
  non-0, and have respective weights $\lambda,\lambda +
  2s,\ldots,\lambda+2sp$. The proposition follows.
\end{proof}

\begin{rem}
  \label{rem:estimates-remark}
  The following analogue of the proposition for $H_n = \XX_n(\Witt_n)
  \cdot T \le \Gamma_n$ may be proved by the same method: if $(\rho,V)$ is
  an $H_n$ module such that $\rho(\XX_n(0,\dots,0,1)) \not = 1$, then
  there are weights $s \in \Z_{>0}$ and $\lambda \in \Z$ such that
  $V_{\lambda + 2sj} \not = 0 $ for $0 \le j \le p^{n-1}$. In
  particular, $T$ has at least $p^{n-1} +1$ distinct weight spaces on
  $V$.
\end{rem}

\section{Minimality of $p+3$}
\label{sec:minimality}

In this section, $G$ again denotes the group $\Gamma_2 = \SL_2(\Witt_2)$,
and $\XX = \XX_2$.

\begin{lem}
  \label{lem:1-dim-ws}
  Let $(\rho,V)$ be a $G$-representation with $\rho(\ZZ(1)) \not =
  1_V$. For some $\nu \in \Z$, the $T$-weight space $V_\nu$ must
  satisfy $\dim V_\nu \ge 2$.
\end{lem}

\begin{proof}
  We may  find $\nu \in \Z$ and a $T$-weight
  vector $v \in V_\nu$ for which
  \begin{equation*}
    \rho(\ZZ(1))v \not = v.
  \end{equation*}  
  There are uniquely determined vectors $v=v_0,v_1,\dots,v_N \in V$
  with $\rho(\ZZ(s))v = \sum_{i=0}^N s^i v_i$ and $v_N \not = 0$.
  Since $\rho(\ZZ(1))v \not = v$, we must have $N>1$. Since
  $\rho(\ZZ(1))v_N = v_N$,  the vectors $v$ and $v_N$ are linearly
  independent.  By \eqref{eq:centralize} we have $v_N \in V_\nu$,
  whence the lemma.
\end{proof}

\begin{theorem}
  \label{theorem:trivial-action}
  Suppose that $(\rho,V)$ is a $G$-representation with $\dim V \le
  p+2$.  Then $\rho(\ZZ(1))=1_V$. In particular, any faithful
  $G$-representation has dimension at least $p+3$.
\end{theorem}

\begin{proof}
  Let $(\rho,V)$ be a $G$-representation for which $\rho(\ZZ(1)) \not
  = 1_V$. By Lemma \ref{lem:generate}, we have $\rho(\XX(0,1)) \not
  = 1_V.$ According to Proposition \ref{prop:weight-estimates} we may
  find $\lambda \in \Z$ and $s > 0$ such that $V_{\lambda + 2sj} \not =
  0$ for $0 \le j \le p$. Since by Lemma \ref{lem:1-dim-ws} there must
  be some $\mu \in \Z$ with $\dim V_\mu \ge 2$, we deduce that $\dim V
  \ge p+2$.
  
  To finish the proof, we suppose that $\dim V = p+2$ and deduce a
  contradiction.  Since we may suppose that $V$ has a 2 dimensional
  weight space $V_\mu$, we see that the $T$-weights of $V$ are
  precisely the $\lambda + 2sj$ for $0 \le j \le p$.  Since the
  character of $V$ must be the character of an $\SL_2(k)$ module, we
  have $\dim V_\gamma = \dim V_{-\gamma}$ for all weights $\gamma \in
  \Z$. Since $V_\mu$ is the unique 2 dimensional weight space, we
  deduce that $\mu = 0$.
  
  It follows that $\lambda,\lambda + 2s,\dots,\lambda + 2sp$ must be the
  weights of some $\SL_2(k)$ module. Steinberg's tensor product
  theorem now implies that $s = p^r$ for some $r \ge 0$. We then have
  $\lambda = -(\lambda + 2p^{r+1}),$ so that $\lambda = -p^{r+1}$.  If
  $p>2$, then we see that $\lambda + 2p^rj \not = 0$ for any $j$, so
  $0$ is not a weight of $V$; this gives our contradiction when $p>2$.
  
  So we may suppose that $p=2$, that $\dim V_{\pm 2^r} = 1$, and that
  $\dim V_0 = 2$. Thus the composition factors of $V$ are $L(2^r) =
  L(1)^{[r]}$, $L(0)$, and $L(0)$. We claim first that $\dim V^R = 1$.
  Indeed, since $\rho(\XX(0,1)) \not = 1$, a look at the proof of
  Theorem \ref{prop:weight-estimates} shows that $V_{\pm2^r} \cap V^R
  = 0$.  Moreover, since $\rho(\ZZ(1)) \not = 1$, Lemma
  \ref{lem:1-dim-ws} shows that $V_0 \not \subset V^R$.
  
  Next, we claim that $\soc(V/V^R)$ can not have $L(0)$ as a summand.
  Indeed, otherwise one finds a 2 dimensional indecomposable
  $G$-module with composition factors $L(0), L(0)$ on which
  $\ZZ(1)$ acts non-trivially. But $\XX(1,0)$ must act trivially
  on such a module, contrary to Lemma \ref{lem:generate}.
  
  It now follows that $\soc(V/V^R) = L(1)^{[r]}$. But then the inverse
  image $W$ in $V$ of $\soc(V/V^R)$ is a $G$-submodule of $V$
  containing $V_{\pm 2^r}$. Moreover, $\dim W = 3$, $\XX(0,1)$ acts
  non-trivially on $W$, while $\ZZ(1)$ must act trivially on $W$.
  Thus $\ker \rho \cap R$ is precisely $Z = \{\ZZ(t) \mid t \in k\}$;
  see Remark \ref{rem:normal-subgroups}.  Let $\YY:\Witt_2 \to \Gamma_2$ be
  the map
  \begin{equation*}
    \YY(w) = \
    \begin{pmatrix}
      1 & 0 \\
      w & 1
    \end{pmatrix}.
  \end{equation*}
  Since $\YY(0,1) \not \in R$, we have $\rho(\YY(0,1)) \not = 1$.
  Moreover, we know that $\rho(\YY(0,1))$ commutes with
  $\rho(\XX(0,1))$. But the fixed point space of $\rho(\XX(0,1))$ on
  $W$ is precisely $W_0 \oplus W_2$, which is not stable under
  $\rho(\YY(0,1))$ by (the proof of) Proposition
  \ref{prop:weight-estimates}. This gives the desired contradiction
  when $p=2$.
\end{proof}

\begin{cor}
  Suppose that $(\rho,V)$ is a $G$-representation with $\dim V \le
  p+2$.  Then the $p$-th power of each unipotent element of $G$ acts
  trivially on $V$.
\end{cor}
\begin{proof}
  Theorem \ref{theorem:trivial-action} shows that $R \cap
  \ker(\rho)$ is a normal subgroup of $G$ containing $\ZZ(1)$, hence
  is $R$ by Lemma \ref{lem:generate}.  If $u \in G$ is unipotent, then
  $u^p \in R$ whence the corollary.
\end{proof}

\section{The Lie algebra of $\Gamma_2$}
\label{sec:inf-faithful}

Let $\glie = \Lie(\Gamma_2)$.  There is an exact sequence of $p$-Lie algebras
and of $\Gamma_2$-modules
\begin{equation}
  \label{eq:Lie-exact}
   0 \to \Lie(R) \to \glie \to \lie{sl}_2(k) \to 0.
\end{equation}

\begin{lem}
  \label{lem:R-acts-trivially}
  Suppose that $p>2$. Then $R$ acts trivially on $\glie$. In particular,
  \eqref{eq:Lie-exact} is an exact sequence of $\SL_2(k)$-modules.
\end{lem}

\begin{proof}
  Since the adjoint module for $\SL_2(k)$ is simple when $p>2$, it
  suffices by Lemma \ref{lem:generate} to show that $\Ad(\XX(0,1))=1$.
  Note that the weights of $T$ on $\glie$ are $\pm 2$, $\pm 2p$, and $0$.
  Since $p>2$, Proposition \ref{prop:weight-estimates} implies that
  $\XX(0,1)$ acts trivially on $\glie$ as desired.
\end{proof}

The Abelian Lie algebra $\Lie(\Witt_2)$ contains an element $Y$ for
which $Y$ and $Y^{[p]}$ form a $k$-basis. The element $Y^{[p]}$ spans
the image of the differential of $(t \mapsto (0,t)):k \to \Witt_2$.
Write $X = d\XX(Y)$. Then $X \not \in \Lie(R)$ and $X^{[p]} \in
\Lie(R)$.

\begin{prop}
  \begin{enumerate}
  \item If $(d\lambda,V)$ is a restricted representation of the
    $p$-Lie algebra $\glie$, then $\ker d\lambda \cap \Lie(R) = 0$ if
    and only if $(d\lambda)(Z) \not = 0$ where $Z=d\ZZ(1)$.
  \item Let $p>2$. Then  \eqref{eq:Lie-exact} is
    split as a sequence of $\Gamma_2$-modules.
  \item Let $p>2$, and let $(d\lambda,V)$ be a representation of
    $\glie$ as a $p$-Lie algebra.  Then $\ker d\lambda = 0$ if and
    only if $d\lambda(X^{[p]}) \not = 0$.
  \end{enumerate}
\end{prop}

\begin{proof}
  (1) is a consequence of Lemma \ref{lem:generate}.
  
  (2) By Lemma \ref{lem:R-acts-trivially}, $R$ acts trivially on
  $\glie$, so $\glie$ may be viewed as a module for $\SL_2(k)$. Note
  that \eqref{eq:Lie-exact} has the form $0 \to L(2p) \to \glie \to
  L(2) \to 0$.  Since $p>2$, $2$ and $2p$ are not linked under the
  action of the affine Weyl group. Hence, the sequence splits thanks
  to the linkage principle \cite[II.6.17]{JRAG}. 
  
  (3) By hypothesis both $d\lambda(X)$ and $d\lambda(X^{[p]})$ are
  non-0.  The image of $X$ is a generator for $\glie/\Lie(R)$ as a
  $\Gamma_2$-module, and $X^{[p]}$ is a generator for $\Lie(R)$ as a
  $\Gamma_2$-module, so the claim follows from (2).
\end{proof}

\begin{cor}
  Consider the $\Gamma_2$ representation $(\rho_p,\A_p)$ of \S \ref{sec:faithful-example}.
  \begin{enumerate}
  \item If $p>2$, then $(d\rho_p,\A_p)$ is a faithful representation
    of $\glie$.
  \item If $p=2$, then $\ker d\rho_2 = \Lie(T)$ is 1-dimensional.
  \end{enumerate}
\end{cor}

\begin{proof}
With notations as before,
using \eqref{eq:witt-vector-ops} one sees that
\begin{equation*}
  \rho_p(\XX(0,s)) A_1 = A_1 + sB_0^p
\end{equation*}
for $s \in k$. It follows that $d\rho_p(X^{[p]}) A_1 = cB_0^p$ for
some $c \in k^\times$. When $p>2$, part (3) of the proposition shows
that $d\rho_p$ is faithful.

Let $Z = d\ZZ(1)$ as before. The calculation in the proof of Theorem
\ref{theorem:faithful} implies that $d\rho_p(Z)A_1 = A_0^p$.  In
particular, part (1) of the proposition shows that $\ker d\rho_p \cap
\Lie(R) = 0$ for all $p$. When $p=2$, note that $\Lie(T)$ indeed acts
trivially; see Remark \ref{rem:faithful-remarks} (a). The corollary
now follows.
\end{proof}

\section{Representations of the associated finite groups}
\label{sec:finite-case}

In this section, a representation of a group is always assumed to be on a finite
dimensional $k$-vector space.

\subsection{}

Let $C$ be a finite cyclic group of order relatively prime to $p$, and
suppose that $\rho:C \to \Aut_{\text{$k$-alg}}(A)$ is a representation
of $C$ by algebra automorphisms on the algebra of truncated
polynomials
\begin{equation*}
  A = k[z]/(z^N)
\end{equation*}
for some $N \ge 2$.  Let $X = \Hom(C,k^\times)$ be the group of
characters of $C$. Since $|C|$ is prime to $p$, $X$ is
(non-canonically) isomorphic to $C$; in particular, it is cyclic.
Note that an element $\mu \in X$ is a generator if and only if $\mu$
is injective as a homomorphism.  
If $(\rho,V)$ is a $C$-representation, and $\mu \in X$, let 
\begin{equation*}
  V_\mu = \{v\in V \mid \rho(c)v = \mu(c)v \ \text{for each $c \in C$}\}.  
\end{equation*}
Of course, $V \iso \bigoplus_{\mu \in X} V_\mu$.

Write $\mm = (z)$ for the maximal
ideal of $A$.

\begin{lem}
  \label{lem:trunc-poly-lemma}
  With notations as above, if $(\rho,A)$ is a faithful
  $C$-representation, then there is $\mu \in X$ and an element $f \in
  \mm \cap A_\mu$  such that $f$ has non-zero image in $\mm/\mm^2$.
\end{lem}

\begin{proof}
  Since $C$ acts by algebra automorphisms, the ideal $\mm^i$ is
  $C$-invariant for each $i \ge 1$.  Since the $C$ representation
  $(\rho,\mm)$ is semisimple, the subrepresentation $\mm^2$ has a
  complement $k\cdot f$ for some $0 \not = f \in \mm$. Thus there is
  $\mu \in X$ such that $\rho(c)f = \mu(c)f$ for each $c \in C$, and
  since $f \not \in \mm^2$, the image of $f$ in $\mm/\mm^2$ is
  non-zero. It remains to argue that $\mu$ is a generator for $X$.
  Note that $1,f,f^2,\ldots,f^{N-1}$ form a $k$-basis for $A$, so that
  \begin{equation*}
    (\rho,A) \iso 1 \oplus \mu \oplus \mu^2 \oplus \cdots \oplus \mu^{N-1}
  \end{equation*}
  as $C$-representations. Since $(\rho,A)$ is a faithful
  representation, we see that $\mu$ must itself be a faithful
  representation of $C$, so that $\mu$ indeed generates $X$.
\end{proof}

\subsection{}
\label{sub:reg-elts} 

Let $V$ be a $k$-vector space of dimension $n \ge 2$.  Let $u$ be a
regular unipotent element in $\GL(V)$; thus $u$ acts on $V$ as a
single unipotent Jordan block. It is well known (and easy to see) that
the centralizer of $u$ in $\lie{gl}(V)= \End_k(V)$ is the
(associative) algebra $k[u]$ generated by $u$.  Let $A = u-1$. Then
$A$ is a regular nilpotent element (it acts as a single nilpotent Jordan
block), and $k[u] = k[A]$. Now, $k[A]$ is isomorphic to the algebra of
truncated polynomials $k[z]/(z^{n-1})$. Moreover, $f \in k[A]$ is a
regular nilpotent element of $\lie{gl}(V)$ if and only if $f \in \mm
\setminus \mm^2.$

\subsection{}

Suppose that $H$ is a finite group, that $C < H$ is a cyclic subgroup
of order prime to $p$, and that $W < H$ is an Abelian $p$-group which
is normalized by $C$. As before, let $X = \Hom(C,k^\times)$. Write $C'$ for
the centralizer in $C$ of $W$, and let
$X' = \{\mu \in X \mid \mu_{\mid C'} = 1\}$.

\begin{prop}
  \label{prop:finite-weight-spaces}
  Let $(\rho,V)$ be a faithful, finite dimensional $H$-representation,
  and suppose that $\rho(W)$ contains a regular unipotent element of
  $\GL(V)$.  If $|C/C'| \ge \dim V$, then $V_\mu$ is 1 dimensional for
  each $\mu \in X$. Moreover, there is $\lambda \in X$ and a generator
  $\mu \in X'$ such that $V \iso V_{\lambda} \oplus V_{\lambda + \mu} \oplus
  \cdots \oplus V_{\lambda + d\mu}$ where $\dim V = d+1$.
\end{prop}

\begin{proof}
  Let $u \in \rho(W)$ be a regular unipotent element. As in
  \ref{sub:reg-elts}, the centralizer of $u$ in $\End_k(V)$ is $k[u]$.
  For each $c \in C$ we have, $\rho(c)u\rho(c)^{-1} \in \rho(W)
  \subset k[u]$ since $W$ is Abelian. It follows that $C$ acts by
  conjugation on $A = k[u]$.  Moreover, $C/C'$ acts faithfully on $A$.
  According to Lemma \ref{lem:trunc-poly-lemma}, there is a generator
  $\mu \in X' = X(C/C')$ and (in view of \ref{sub:reg-elts}) a regular
  nilpotent element $A \in (\lie{gl}(V))_\mu$.
  
  Let $d = \dim V -1$. We may thus find $\lambda \in X$ such that
  $A^d(V_\lambda) \not = 0$.  It follows that $V_\lambda, V_{\lambda +
    \mu},\ldots, V_{\lambda + d\mu}$ are all non-0.  Since $\mu$ has
  order $|C/C'| > d$, each of these subspaces has dimension 1. The
  proposition follows.
\end{proof}

\subsection{}

Fix a $p$-power $q = p^a$, and let $\F_q$ be the field with $q$
elements.  The group $\Gamma_n$, and the homomorphisms $\phi:\G_m \to
\Gamma_n$ and $\XX_n:\Witt_n \to \Gamma_n$ are defined over $\F_q$.

Let $n=2$, and let $C,W \le \Gamma_2(\F_q) = \SL_2(\Witt_2(\F_q))$ be
respectively the image under $\phi$ of $\G_m(\F_q) \iso
\F_q^\times$, and  the image under $\XX_2$ of $\Witt_2(\F_q)$.
Then $C$ is cyclic of order prime to $p$, and $W$ is a $p$-group normalized
by $C$. Moreover, the centralizer $C'$ of $W$ in $C$ has order 2.

\begin{theorem}
  Suppose that $p \ge 3$ and $q \ge p^2$.  Then the minimal dimension
  of a faithful $k$-representation of $\Gamma_2(\F_q)$ is $p+3$.
\end{theorem}

\begin{proof}
  That $\Gamma=\Gamma_2(\F_q)$ has a faithful representation of dimension $p+3$ follows
from Remark \ref{rem:faithful-remarks}(b). 

We now suppose that $(\rho,V)$ is a faithful representation of
$\Gamma$ with $\dim V \le p+2$ and deduce a contradiction.  Since the
element $\XX_2(1,0)$ of $\Gamma_2(\F_q)$ has order $p^2$, we see that
$\dim V \ge p+1$. Suppose first that $\dim V = p+1$.  Then the image
$\rho(W)$ must contain a regular unipotent element.  

With our assumption on $q$, we have $|C/C'| = \dfrac{q-1}{2} \ge p+1 =
\dim V$.  An application of Proposition
\ref{prop:finite-weight-spaces} for the subgroups $C,W < \Gamma$
therefore shows that the spaces $V_\mu$ with $\mu \in X(C)$ are all
1-dimensional. Lemma \ref{lem:1-dim-ws} now shows that the element
$\ZZ(1) \in \Gamma$ must act trivially; this contradicts our
assumption that $(\rho,V)$ is faithful.

Finally, suppose that $\dim V = p+2$. Let $H$ be the subgroup of
$\Gamma$ generated by $C$ and $W$. Since $H$ is nilpotent and since
$\rho(H)$ contains a unipotent element with Jordan block sizes
$(p+1,1)$, we have $V = V' \oplus V''$ with $V'$ and $V''$ invariant
under $H$, and with $\dim V' = p+1$.  Now an application of
Proposition \ref{prop:finite-weight-spaces} to the $H$-representation
$V'$ shows that there is a weight $\lambda \in X(C)$ and a generator
$\mu \in X'$ such that $V' = \bigoplus_{i=0}^p V'_{\lambda +
  i\mu}$ with $\dim V'_{\lambda + i \mu} =1$ for each $i$. In view of
Lemma \ref{lem:1-dim-ws}, there is precisely one $\gamma \in X$
with $\dim V_\gamma = 2$.

As in \S \ref{sec:unip-radical}, the composition factors of the
$\Gamma$-representation $(\rho,V)$ may be identified with simple
representations of the group $\SL_2(\F_q) = \Gamma/R(\F_q)$.  Thanks
to a theorem of Curtis (\cite[Theorem 43]{Steinberg}) the
semisimplification of $(\rho,V)$ is the restriction to $\SL_2(\F_q)$
of a semisimple rational $\SL_2(k)$ module $(\psi,W)$ with $\dim W =
p+2$ and with precisely one two-dimensional weight space.  As in the
proof of Theorem \ref{theorem:trivial-action}, one knows that this is
impossible (since $p > 2$).
\end{proof}

\begin{example}
  As a ``concrete'' example, let $A = \Z[i]$ be the ring of Gaussian
  integers. Suppose that the prime $p$ satisfies $p \congruent 1 \pmod
  4$; such a prime may be written $p = a^2 + b^2$ for $a,b\in \Z$.
  Denoting by $\mathfrak{P}$ the ideal $(a + bi)A$, one has
  $A/\mathfrak{P} \iso \F_{p^2}$. Then $A/\mathfrak{P}^2 \iso
  \Witt_2(\F_{p^2})$, so the minimal dimension of a faithful
  $p$-modular representation of $\SL_2(A/\mathfrak{P^2})$ has
  dimension $p+3$.
\end{example}

\providecommand{\bysame}{\leavevmode\hbox to3em{\hrulefill}\thinspace}

\end{document}